\documentclass[twoside,11pt,reqno]{article}
\usepackage{amssymb}
\usepackage{amsthm}
\usepackage[tbtags]{amsmath}
\usepackage{doc}
\usepackage{latexsym}
\usepackage{amscd}
 \font\Bbb=msbm10 \theoremstyle{change}

\theoremstyle{change}  
\begin{document} \thispagestyle{plain}
 \markboth{}{}
\small{\addtocounter{page}{0} \pagestyle{plain} \vspace{0.2in}
\noindent{\large \bf A New Changhee $q$-Euler Numbers and
Polynomials Associated with $p$-Adic $q$-Integrals}
\footnote{{}
\indent 2000 Mathematics Subject Classification: $p$-adic
$q$-integral, $q$-Euler numbers and polynomials.}
\vspace{0.15in}\\
\noindent{\sc Seog-Hoon Rim}
\newline
{\it Department of Mathematics Education, Kyungpook National University, Daegu 702-701, Korea\\
e-mail} : {\verb|shrim@knu.ac.kr|}
\vspace{0.15in}\\
\noindent{\sc Taekyun Kim}
\newline
{\it EECS, Kyungpook National University, Taegu 702-701, S.Korea \\
e-mail} : {\verb|tkim@knu.ac.kr, tkim64@hanmail.net|}
\vspace{0.15in}\\

{\footnotesize {\sc Abstract.}  Using non-archimedean $q$-integrals
on $\mathbb{Z}_{p}$ defined in [15, 16], we define a new Changhee
$q$-Euler polynomials and numbers which are different from those of
Kim [7] and Carlitz [2]. We define generating functions of multiple
$q$-Euler numbers and polynomials. Furthermore we construct
multivariate Hurwitz type zeta function which interpolates the
multivariate $q$-Euler numbers or polynomials at negative integers.}
\vspace{0.2in}\\
\noindent{\bf 1. Introduction} \setcounter{equation}{0}
\vspace{0.1in}\\
\indent Let $p$ be a fixed odd prime in this paper. Throughout
this paper, the symbols $\mathbb{Z}$, $\mathbb{Z}_{p}$,
$\mathbb{Q}_{p}$, $\mathbb{C}$ and $\mathbb{C}_{p}$, denote the
ring of rational integers, the ring of $p$-adic integers, the
field of $p$-adic numbers, the complex number field, and the
completion of the algebraic closure of $\mathbb{Q}_{p}$,
respectively. Let $\nu_{p}(p)$ be the normalized exponential
valuation of $\mathbb{C}_{p}$ with $|p|_{p}=p^{-\nu_{p}(p)}=
p^{-1}$. When one speaks of $q$-extension, $q$ can be regarded as
an indeterminate, a complex number $q \in \mathbb{C}$, or a
$p$-adic number $q \in \mathbb{C}_{p}$; it is always clear from
the context. If  $q \in \mathbb{C}$, then one usually assumes that
$|q|<1$. If $q \in \mathbb{C}_{p}$, then one usually assumes that
$|q-1|_{p} < p^{-\frac{1}{p-1}}$, and hence $q^{x} = \exp ( x \log
q)$ for $x \in \mathbb{Z}_{p}$. In this paper, we use the below
notation
$$[x]_{q} = \dfrac{1-q^{x}}{1-q}, \quad (a:q)_{n} =
(1-a)(1-aq)\cdots(1-aq^{n-1}), \quad {\rm cf.~[3,4,5,6,7,8]}.$$
Note that $\lim_{q\rightarrow 1}[x]_{q} = x$ for any $x$ with
$|x|_{p} \leq 1$ in the $p$-adic case. For a fixed positive
integer $d$ with $(p, d)=1$, set
\begin{eqnarray*}
&&X=X_{d}=\lim_{\substack{\longleftarrow \\ N}} \mathbb{Z}/dp^{N},\\
&&X_{1} = \mathbb{Z}_{p}, \quad X^{*}=
\bigcup_{\substack{0<a<dp \\ (a, p)=1}} a+dp\mathbb{Z}_{p},\\
&&a+dp^{N}\mathbb{Z}_{p} = \{x \in X ~|~ x \equiv a ~({\rm
mod}~p^{N})\},
\end{eqnarray*}
where $a\in \mathbb{Z}$ satisfies the condition $0 \leq a <
dp^{N}$, (cf. [10,11]). We say that $f$ is a uniformly
differentiable function at a point $a \in \mathbb{Z}_{p}$, and
write $f \in UD(\mathbb{Z}_{p})$, if the difference quotients
$F_{f}(x, y) = \frac{f(x)-f(y)}{x-y}$ have a limit $f'(a)$ as $(x,
y) \rightarrow (a, a)$. For $f \in UD(\mathbb{Z}_{p})$, let us
begin with the expression
$$\dfrac{1}{[p^{N}]_{q}} \sum_{0 \leq j < p^{N}} q^{j} f(j) = \sum_{0 \leq j <
p^{N}}f(j)\mu_{q}(j+p^{N}\mathbb{Z}_{p}), \quad {\rm
cf.~[4,5,6,8]},$$ which represents a $q$-analogue of Riemann sums
for $f$. The integral of $f$ on $\mathbb{Z}_{p}$ is defined as the
limit of those sums(as $n\rightarrow \infty$) if this limit
exists. The $p$-adic $q$-integral of a function $f \in
UD(\mathbb{Z}_{p})$ is defined by
$$I_{q}(f) = \int_{X} f(s) d\mu_{q} (x) = \int_{X_{d}} f(x)
d\mu_{q}(x) = \lim_{N \rightarrow \infty}\dfrac{1}{[dp^{N}]_{q}}
\sum_{x=0}^{dp^{N}-1} f(x) q^{x}.$$
\vspace{0.02in}\\
Recently, many mathematicians studied Bernoulli and Euler numbers
(see [1-24]). Using non-archimedean $q$-integrals on
$\mathbb{Z}_{p}$ defined in [15, 16], we define a new Changhee
$q$-Euler polynomials and numbers which are different from those of
Kim [7] and Carlitz [2]. We define generating functions of multiple
$q$-Euler numbers and polynomials. Furthermore we construct
multivariate Hurwitz type zeta function which interpolates the
multivariate $q$-Euler numbers or polynomials at negative integers.

\noindent{\bf 2. Multivariate $q$-Euler numbers and polynomials}
\vspace{0.1in}\\
\indent Using $p$-adic $q$-integrals on $\mathbb{Z}_{p}$, we now
define a new $q$-Euler polynomials as follows:
$$\int_{\mathbb{Z}_{p}} q^{y}(x+y)^{n} d\mu_{-1}(y) = E_{n,
q}(x).$$ Note that $\lim_{q\rightarrow 1} E_{n, q}(x) = E_{n}(x)$,
where $E_{n}(x)$ are Euler polynomial which are defined by
$\frac{2}{e^{t}+1} e^{xt} = \sum_{n=0}^{\infty} E_{n}(x)
\frac{t^{n}}{n!}$. In the case $x=0$, $E_{n, q}=E_{n, q}(0)$ are
called a new $q$-Euler numbers. And note that $\lim_{q\rightarrow
1}E_{n, q}=E_{n}$, where $E_{n}$ are classical Euler numbers. Let
$a_{1}, a_{2}, \cdots, a_{r}, b_{1}, b_{2}, \cdots, b_{r}$ be
positive integers. Then we consider a multivariate integral as
follows
\begin{eqnarray}
&&\int_{\mathbb{Z}_{p}}\cdots
\int_{\mathbb{Z}_{p}} q^{b_{1}x_{1} + \cdots + b_{r}x_{r}}
(a_{1}x_{1} + \cdots + a_{r}x_{r})^{n} d\mu_{-1}(x_{1})\cdots
d\mu_{-1}(x_{r})\notag\\
&=& {E_{n, q}}^{(r)} (x| a_{1}, a_{2}, \cdots, a_{r} ~;~ b_{1},
b_{2}, \cdots, b_{r}).
\end{eqnarray}
Here ${E_{n, q}}^{(r)} (x| a_{1}, a_{2}, \cdots, a_{r} ~;~ b_{1},
b_{2}, \cdots, b_{r})$ are called multivariate $q$-Euler
polynomial of order $r$.

\smallskip

In the case $x=0$ in (1), $E_{n, q} (0| a_{1}, a_{2}, \cdots,
a_{r} ~;~ b_{1}, b_{2}, \cdots, b_{r})= E_{n, q}(a_{1}, a_{2},
\cdots, a_{r} $ $~;~ b_{1}, b_{2}, \cdots, b_{r})$ will be called
multivariate $q$-Euler numbers of order $r$.

From (1), we derive the following generating function for
multivariate $q$-Euler polynomials
\begin{eqnarray}
&&F_{q}^{(r)}(t, x| a_{1}, a_{2}, \cdots, a_{r} ~;~ b_{1}, b_{2},
\cdots, b_{r})\notag\\
&=& \sum_{n=0}^{\infty} {E_{n, q}}^{(r)} (x| a_{1}, a_{2}, \cdots,
a_{r} ~;~ b_{1}, b_{2}, \cdots, b_{r}) \dfrac{t^{n}}{n!}\notag\\
&=& \dfrac{2^{r}}{(q^{b_{1}}e^{a_{1}t}+1)
(q^{b_{2}}e^{a_{2}t}+1)\cdots (q^{b_{r}}e^{a_{r}t}+1)}e^{xt}.
\end{eqnarray}
Note that ${E_{0, q}}^{(r)} (x| a_{1}, a_{2}, \cdots, a_{r} ~;~
b_{1}, b_{2}, \cdots, b_{r})=
\dfrac{2^{r}}{[2]_{q^{b_{1}}}[2]_{q^{b_{2}}} \cdots
[2]_{q^{b_{r}}}}.$

\medskip

Let $\chi$ be the Dirichlet character with conductor $f$(=odd)$\in
\mathbb{N}$. Then we define the multivariate generalized $q$-Euler
numbers attached to $\chi$ as follows:
\begin{eqnarray}
&&\int_{X}\cdots \int_{X} \chi(x_1)\cdots \chi(x_r) q^{b_{1}x_{1} +
\cdots + b_{r}x_{r}} (a_{1}x_{1} + \cdots + a_{r}x_{r})^{n}
d\mu_{-1}(x_{1})\cdots d\mu_{-1}(x_{r})\notag\\
&=& E_{n, \chi, q} ( a_{1}, \cdots, a_{r} ~;~ b_{1}, \cdots,
b_{r}).
\end{eqnarray}
From (3), we can derive generating function for the multivariate
generalized $q$-Euler numbers attached to $\chi$,
\begin{eqnarray}
&&F_{\chi, q}^{(r)}(t | a_{1}, \cdots, a_{r} ~;~ b_{1}, \cdots,
b_{r})\notag\\
&=&\sum_{n=0}^{\infty} E_{n, \chi, q} (a_{1}, \cdots,
a_{r} ~;~ b_{1}, \cdots, b_{r}) \dfrac{t^{n}}{n!}\notag\\
&=& \sum_{n_{1}, \cdots, n_{r}=0}^{f-1}
\dfrac{2^{r}(-1)^{n_{1}+\cdots + n_{r}}q^{b_{1}n_{1}+\cdots +
b_{r}n_{r}} e^{(a_{1}n_{1}+\cdots + a_{r}n_{r})t}
\chi(n_{1})\cdots \chi(n_{r})}{(q^{fb_{1}}e^{fa_{1}t}+1)\cdots
(q^{fb_{r}}e^{fa_{r}t}+1)},
\end{eqnarray}
where $n$ is odd.

\medskip

Let $n$ be odd positive integer. Form the definition of ${E_{n,
\chi, q}}^{(r)}(a_{1}, \cdots, a_{r} ~;~ b_{1}, \cdots, b_{r})$ in
(3), we have the following;
\begin{eqnarray*}
&&{E_{n, \chi, q}}^{(r)}(a_{1}, \cdots, a_{r} ~;~ b_{1}, \cdots,
b_{r})\\
&=&\lim_{N \rightarrow \infty} \sum_{n_{1}, \cdots,
n_{r}=0}^{fp^{N}-1} \chi(n_{1})\cdots \chi(n_{r})
(-1)^{n_{1}+\cdots+n_{r}} q^{b_{1}n_{1}+\cdots+b_{r}n_{r}}(a_{1}x_{1}+\cdots + a_{r}x_{r})^{n}\\
&=&\lim_{N \rightarrow \infty} \sum_{n_{1}, \cdots,
n_{r}=0}^{f-1}\sum_{x_{1}, \cdots, x_{r}=0}^{p^{N}-1}
\Pi_{j=1}^{r} \chi(n_{j}+fn_{j}) (-1)^{\sum_{j=1}^{r}
(n_{j}+fx_{j})}q^{\sum_{j=1}^{r} b_{j}(n_{j}+x_{j}f_{1})}\\
&&~~ \sum_{j=1}^{r} a_{j} (n_{j}+fx_{j})^{n} \\
&=&f^{n}
\sum_{n_{1}, \cdots, n_{r}=0}^{f-1} (-1)^{\sum_{j=1}^{r}
n_{j}} \Pi_{j=1}^{r} \chi(n_{j}) q^{\sum_{j=1}^{r} b_{j}n_{j}}\\
&&~\lim_{N \rightarrow \infty} \sum_{x_{1}, \cdots,
x_{r}=0}^{p^{N}-1} (-1)^{\sum_{j=1}^{r} (x_{j})} q^{f
\sum_{j=1}^{r}b_{j}x_{j}} \left(\dfrac{\sum_{j=1}^{r}
a_{j}n_{j}}{f} + \sum_{j=1}^{r}
a_{j}x_{j}\right)^{n}
\end{eqnarray*}
\begin{eqnarray*}
&=&f^{n} \sum_{n_{1}, \cdots, n_{r}=0}^{f-1} (-1)^{\sum_{j=1}^{r}
n_{j}} \Pi_{j=1}^{r} \chi(n_{j}) q^{\sum_{j=1}^{r} b_{j}n_{j}}\\
&&~~\underbrace{\int_{\mathbb{Z}_{r}}\cdots
\int_{\mathbb{Z}_{p}}}_{r-{\rm times}} q^{f \sum_{j=1}^{r}
b_{j}x_{j}} \left(\dfrac{\sum_{j=1}^{n} a_{j}x_{j}}{f} +
\sum_{j=1}^{n} a_{j}x_{j}\right)^{n} d\mu_{-1}(x_{1}) \cdots
d\mu_{-1}(x_{r})\\
&=&f^{n} \sum_{n_{1}, \cdots, n_{r}=0}^{f-1}
(-1)^{\sum_{j=1}^{r}
n_{j}} \Pi_{j=1}^{r} \chi(n_{j}) q^{\sum_{j=1}^{r} b_{j}n_{j}}\\
&&~~E_{n, q^{f}} \left(\dfrac{\sum_{j=1}^{n} a_{j}x_{j}}{f}
~\Big|~ a_{1}, \cdots, a_{r} ~;~ b_{1}, \cdots, b_{r} \right).
\end{eqnarray*}

Therefore we obtain;
\vspace{0.1in}\\
\noindent {\bf Theorem 1.} {\it Let $a_{1}, \cdots, a_{r}, b_{1},
\cdots, b_{r}$ be positive integers then we have
\begin{eqnarray}
&&{E_{n, \chi, q}}^{(r)}(a_{1}, \cdots, a_{r} ~;~ b_{1}, \cdots,
b_{r})\notag\\
&=& \sum_{n_{1}, \cdots, n_{r}=0}^{f-1} (-1)^{\sum_{j=1}^{r}
n_{j}} \Pi_{j=1}^{r} \chi(j) q^{\sum_{j=1}^{r} b_{j}n_{j}}E_{n,
q^{f}} (a_{1}, \cdots, a_{r} ~;~ b_{1}, \cdots,
b_{r}),\end{eqnarray} where $f$, $n$ are odd positive integers.}
\vspace{0.2in}\\
\noindent{\bf 3. Multivariate $q$-zeta functions in $\mathbb{C}$}
\vspace{0.1in}\\
\indent In this section we assume $q \in \mathbb{C}$ with $|q|<1$.
Let us assume that $a_{1}, \cdots, a_{r}, ~ b_{1}, \cdots, b_{r}$
are positive integers.

The purpose of this section is to study multivariate Hurwitz type
zeta function which interpolates multivariate $q$-Euler
polynomials of order $r$ at negative integers. By (2), we easily
see that
\begin{eqnarray}
&&F^{(r)}_{q} (t, x ~|~ a_{1}, \cdots, a_{r} ~;~ b_{1}, \cdots,
b_{r})\notag\\
&=& 2^{r} \sum_{n_{1}, \cdots, n_{r}=0}^{\infty}
(-1)^{\sum_{j=1}^{r} n_{j}} q^{\sum_{j=1}^{r} b_{j}n_{j}}
e^{(\sum_{j=1}^{r} a_{j}n_{j}+x)t}.
\end{eqnarray}
By taking derivatives of order $k$, on both sides of (6) we obtain
the following;
\vspace{0.1in}\\
\noindent {\bf Theorem 2.} {\it Let $k$ be positive odd integer
and let $a_{1}, \cdots, a_{r}, ~b_{1}, \cdots, b_{r}$ be the
positive integers. Then we have}
\begin{eqnarray}
&&{E_{k,q}}^{(r)} (x ~|~a_{1}, \cdots, a_{r}~;~b_{1}, \cdots,
b_{r})\notag\\
&=& 2^{r} \sum_{n_{1}, \cdots, n_{r}=0}^{\infty}
(-1)^{\sum_{j=1}^{r} n_{j}} q^{\sum_{j=1}^{r} b_{j}n_{j}}
(\sum_{j=1}^{r} a_{j}n_{j}+x)^{k}.
\end{eqnarray}
By the above Theorem, we may now construct the complex
multivariate $q$-zeta functions as follows:
\vspace{0.1in}\\
\noindent {\bf Definition 1.} For $s \in \mathbb{C}$, we define
\begin{equation}
\zeta_{r} (s, x ~|~ a_{1}, \cdots, a_{r} ~;~ b_{1}, \cdots, b_{r})
= \sum_{n_{1}, \cdots, n_{r}=0}^{\infty}
\dfrac{2^{r}(-1)^{\sum_{j=1}^{r} n_{j}} q^{\sum_{j=1}^{r}
b_{j}n_{j}}} {(\sum_{j=1}^{r} a_{j}n_{j} +x)^{s}}.
\end{equation}
Thus we note that this function in (8) is analytic continuation in
whole complex plane. And we see this multivariate $q$-zeta
function interpolates $q$-Euler polynomials at negative integers.
\vspace{0.1in}\\
\noindent {\bf Theorem 3.} {\it Let $n$ be an odd positive
integer. Then we have}
\begin{equation}
\zeta_{r} (-n,x ~|~a_{1},\cdots, a_{r}~;~b_{1}, \cdots,
b_{r})={E_{n,q}}^{(r)} (x~|~a_{1}, \cdots, a_{r}~;~b_{1}, \cdots,
b_{r}).
\end{equation}
We now give the complex integral representation of $\zeta_{r} (s,x
~|~a_{1},\cdots, a_{r}~;~b_{1}, \cdots, b_{r})$. Using (6), we
have the following,
\begin{eqnarray}
&&\dfrac{1}{\Gamma(s)} \oint_{\mathbb{C}} F_{q}^{(r)} (-t,x~|~
a_{1}, \cdots, a_{r}~;~b_{1}, \cdots, b_{r}) t^{s-1}dt\notag\\
&=&2^{r} \sum_{n_{1}, \cdots, n_{r}=0}^{\infty}
(-1)^{\sum_{j=1}^{r} n_{j}} q^{\sum_{j=1}^{r} b_{j}n_{j}}
\dfrac{1}{\Gamma(s)} \oint_{\mathbb{C}} e^{-(\sum_{j=1}^{r}
n_{j}a_{j}+x)t} t^{s-1} dt\notag\\
&=&2^{r} \sum_{n_{1}, \cdots, n_{r}=0}^{\infty}
(-1)^{\sum_{j=1}^{r} n_{j}} \dfrac{q^{\sum_{j=1}^{r}
b_{j}n_{j}}}{(\sum_{j=1}^{r} n_{j}a_{j}+x)^{s}}
\dfrac{1}{\Gamma(s)}\oint_{\mathbb{C}} e^{-y} y^{s-1} dy\notag\\
&=& \zeta_{r}(s,x~|~a_{1}, \cdots, a_{r} ~;~b_{1}, \cdots, b_{r}).
\end{eqnarray}
On the other hand
\begin{eqnarray}
&&\dfrac{1}{\Gamma(s)} \oint_{\mathbb{C}} F^{(r)}_{q} (-t,x
~|~a_{1},\cdots, a_{r}~;~b_{1}, \cdots, b_{r}) t^{s-1} dt\notag\\
&=& \sum_{m=0}^{\infty} (-1)^{m}
\dfrac{{E_{m,q}}^{(r)}(x~|~a_{1},\cdots, a_{r}~;~b_{1}, \cdots,
b_{r})}{m!} \dfrac{1}{\Gamma(s)} \oint_{\mathbb{C}} t^{m+s-1} dt.
\end{eqnarray}
Thus by (10) and (11), we have the following:
\begin{eqnarray*}
&&\zeta_{r}(s,x~|~a_{1}, \cdots, a_{r}~;~b_{1}, \cdots,
b_{r})\notag\\
&=& \sum_{m=0}^{\infty} (-1)^{m} \dfrac{E_{m}^{(r)} (x~|~a_{1},
\cdots, a_{r}~;~b_{1}, \cdots,~b_{r})}{m!} \dfrac{1}{\Gamma(s)}
\oint_{\mathbb{C} } t^{m+s-1}dt.
\end{eqnarray*}
Thus we have
\begin{equation}
\zeta_{r} (-n,x~|~a_{1}, \cdots, a_{r}~;~b_{1}, \cdots,
b_{r})=E_{n}^{(r)} (x~|~ a_{1}, \cdots, a_{r}~;~b_{1},
\cdots,b_{r}).
\end{equation}
To construct the multivariate Dirichlet $L$-function we
investigate the generating function of generalized multivariate
$q$-Euler numbers attached to $\chi$, which is derived in (4).
\begin{eqnarray*}
&&2^{r} \sum_{n_{1}, \cdots, n_{r}=0}^{f-1}
\dfrac{(-1)^{\sum_{j=1}^{r} n_{j}} q^{\sum_{j=1}^{r} b_{j}n_{j}}
e^{\sum_{j=1}^{r} a_{j}n_{j}t} \prod_{j=1}^{r}
\chi(n_{j})}{\prod_{j=1}^{r} (q^{fbj} e^{f a_{j} t}+1)}\\
&=& 2^{r} \sum_{n_{1}, \cdots, n_{r}=0}^{f-1} (-1)^{\sum_{j=1}^{r}
n_{j}} q^{\sum_{j=1}^{r} b_{j}n_{j}} \prod_{j=1}^{r} \chi(n_{j})\\
&&~\sum_{x_{1}, \cdots, x_{r}=0}^{\infty} (-1)^{\sum_{j=1}^{r}
n_{j}}
q^{\sum_{j=1}^{r} fb_{j} x_{j}} e^{\sum_{j=1}^{r} f a_{j} x_{j}t}\\
&=&2^{r} \sum_{x_{1}, \cdots, x_{r}=0}^{\infty} \sum_{n_{1},
\cdots, n_{r}=0}^{f-1} (-1)^{\sum_{j=1}^{r} (n_{j}+f x_{j})}
q^{\sum_{j=1}^{r} b_{j}(n_{j} +f x_{j})} \\
&&~\prod_{j=1}^{r}
\chi(n_{j}+f x_{j}) e^{t \sum_{j=1}^{r} a_{j}(n_{j}+f x_{j})}\\
&=&2^{r} \sum_{n_{1}, \cdots, n_{r}=0}^{\infty}
(-1)^{\sum_{j=1}^{r} n_{j}} q^{\sum_{j=1}^{r} b_{j} n_{j}}
\prod_{j=1}^{r} \chi(n_{j}) e^{t \sum_{j=1}^{r} a_{j} n_{j}}.
\end{eqnarray*}
Thus we can write,
\begin{eqnarray}
&&{F_{\chi,q}}^{(r)} (t~|~a_{1}, \cdots, a_{r}~;~b_{1}, \cdots,
b_{r})\notag\\
&=& 2^{r} \sum_{n_{1}, \cdots, n_{r}=1}^{\infty}
(-1)^{\sum_{j=1}^{r} n_{j}} q^{\sum_{j=1}^{r} b_{j} n_{j}}
\prod_{j=1}^{r} \chi (n_{j}) e^{t\sum_{j=1}^{r} a_{j} n_{j}}\notag\\
&&\quad ({\rm since}~ \chi(0) = 0, ~{\rm we ~can ~start ~at}~
n_{1}=1, n_{2}=1, \cdots, n_{r}=1)\notag\\
&=& \sum_{n=1}^{\infty} E_{n,\chi,q}(a_{1}, \cdots, a_{r}~;~b_{1},
\cdots, b_{r}) \dfrac{t^{n}}{n!}. \end{eqnarray}

From (13), we can derive the following
\begin{eqnarray}
&&2^{r} \sum_{n_{1}, \cdots, n_{r}=1}^{\infty}
(-1)^{\sum_{j=1}^{r} n_{j}} q^{\sum_{j=1}^{r} b_{j} n_{j}}
\prod_{j=1}^{r} \chi(n_{j}) (a_{1}n_{1}+ \cdots + a_{r}n_{r})\notag\\
&=&E_{k,\chi,q}(a_{1}, \cdots, a_{r}~;~b_{1}, \cdots, b_{r}).
\end{eqnarray}
Therefore we have the following:
\vspace{0.1in}\\
\noindent {\bf Definition 2.} For $s \in \mathbb{C}$, define
multivariate Dirichlet $L$-function as follows:
\begin{eqnarray}
&&L_{r}(s,x~|~a_{1}, \cdots, a_{r}~;~b_{1}, \cdots, b_{r})\notag\\
&=& 2^{r} \sum_{n_{1}, \cdots, n_{r}=1}^{\infty}
\dfrac{(-1)^{n_{1}+ \cdots + n_{r}} q^{b_{1}n_{2}+\cdots +
b_{r}n_{r}} \chi(n_{1})\cdots \chi(n_{r})}{(a_{1}n_{1}+ \cdots
+a_{r}n_{r})^{s}}.
\end{eqnarray}

Note that $L_{r} (s, \chi~|~ a_{1}, \cdots, a_{r} ~;~ b_{1},
\cdots, b_{r})$ is also analytic function in whole complex plane.
By (12), (13), (14) and (15), we see that the $q$-analogue
multivariate Dirichlet $L$-function interpolates multivariate
generalized $q$-Euler numbers attached to $\chi$ at negative
integers as follows;
\vspace{0.1in}\\
\noindent {\bf Theorem 4.} {\it Let $k$ be a positive integer.
Then we have} $$L_{r} (-k,\chi~|~a_{1}, \cdots, a_{r}~;~b_{1},
\cdots, b_{r}) = E_{k,\chi,q} (a_{1}, \cdots, a_{r}~;~b_{1},
\cdots, b_{r}).$$
\vspace{0.2in}\\
\footnotesize{

\end{document}
\begin{thebibliography}{99}
\bibitem {} L. Carlitz, {\em $q$-Bernoulli numbers and polynomials}, Duke Math. J., {\bf 15}(1948), 987-1000.
\bibitem {} L. Carlitz, {\em $q$-Bernoulli and Eulerian numbers}, Trans. Amer. Math. Soc., {\bf 76}(1954), 332-350
\bibitem {} M. Cenkci, M.Can, {\em Some results on $q$-analogue of the Lerch zeta function},
Adv. Stud. Contemp. Math., {\bf 12}(2006), 213-223
\bibitem {} M. Cenkci, M.Can, V. Kurt {\em $p$-adic interpolation functions and
Kummer-type congruences for $q$-twisted and $q$-generalized twisted
Euler numbers}, Adv. Stud. Contemp. Math., {\bf 9}(2004), 203-216
\bibitem {} H.S. Cho, E. S. Kim, {\em Translation-invariant $p$-adic integral on $\Bbb Z\sb p$},
Notes Number Theory Discrete Math., {\bf 7}(2001), 70-77
\bibitem {} A. S. Hegazi, M. Mansour, {\em A note on $q$-Bernoulli numbers and polynomials
}, J. Nonlinear Math. Phys., {\bf 13}(2006), 9-18
\bibitem {} T. Kim, {\em  $p$-adic $q$-integrals associated with the Changhee-Barnes' $q$-Bernoulli polynomials},
Integral Transforms Spec. Funct., {\bf 15}(2004), 415-420
\bibitem {} T. Kim, S. H. Rim, {\em  Generalized Carlitz's $q$-Bernoulli numbers in the $p$-adic number field},
Adv. Stud. Contemp. Math., {\bf 2}(2000), 9-19
\bibitem {} T. Kim, {\em  On the analogs of Euler numbers and polynomials associated with $p$-adic $q$-integrals on
$\Bbb Z_p$ at $q=-1$}, J. Math. Anal. Appl.(2006)
doi:10.1016/j.jmaa.2006.09.027
 \bibitem {} T. Kim, {\em A new approach to $q$-zeta function}, Adv. Stud. Contemp. Math., {\bf
11}(2005), 157-162
\bibitem {} T. Kim, L. C. Jang, H. K. Park, {\em A note on $q$-Euler and
Genocchi numbers}, Proc. Japan Academy, {\bf 77}(2001), 139-141
\bibitem {} T. Kim, {\em Power series and asymptotic series
associated with the $q$-analog of the two-variable $p$-adic
$L$-function}, Russ. J. Math. Phys., {\bf 12}(2005), 189-196.
\bibitem {} T. Kim, {\em Non-Archimedean $q$-integrals associated
with multiple Changhee $q$-Bernoulli polynomials}, Russ. J. Math.
Phys., {\bf 10}(2003), 91-98.
\bibitem {} T. Kim, {\em An invariant $p$-adic integral
associated with Daehee numbers}, Integral Transforms and special
functions, {\bf 13}(2002), 65-69.
\bibitem {} T. Kim, {\em On a $q$-analogue of the $p$-adic
log gamma functions and relate integrals}, J. Number Theory, {\bf
76}(1999), 320-329.
\bibitem {} T. Kim, {\em $q$-Volkenborn
integration}, Russ. J. Math. Phys., {\bf 9}(2002), 288-299.
\bibitem {}  T. Kim, {\em $q$-Euler numbers and polynomials associated with $p$-adic
$q$-integrals}, J. Nonlinear Math. Phys.(accepted).
\bibitem {} T. Kim, {\em q-Euler numbers and polynomials associated with p-adic q-integrals and basic q-zeta function},
Trends in Mathematics (Information Center for Mathematical
Sciences), {\bf 9}(2006), 7-12.
\bibitem {} B. A. Kupershmidt, {\em Reflection symmetries of $q$-Bernoulli polynomials},
J. Nonlinear Math. Phys., {\bf 12}(2005), 412-422.
\bibitem {} J. Satoh, {\em $q$-analogue of
Riemann's $\zeta$-function and $q$-Euler numbers}, J. Number Theory,
{\bf 31}(1989), 346-362.
\bibitem {} K. Shiratani, S. Yamamoto, {\em On a $p$-adic interpolation function for the Euler numbers
 and its derivatives}, Mem. Fac. Sci. Kyushu Univ. Ser. A,
{\bf 39}(1985), 113-125.
\bibitem {} Y. Simsek, {\em Theorems on twisted $L$-function and twisted Bernoulli numbers},
Adv. Stud. Contemp. Math., {\bf 11}(2005), 205-218.
\bibitem {} Y. Simsek, D. Kim, S. H. Rim, {\em On the two-variable Dirichlet $q$-$L$-series},
Adv. Stud. Contemp. Math., {\bf 10}(2005), 131-142.

\bibitem {} Q.-M. Luo, H. M. Srivastava, {\em Some relationships between the Apostol-Bernoulli
and Apostol-Euler polynomials}, Comput. Math. Appl., {\bf 10}(2005),
631-642.


\end{thebibliography}
